\documentclass[11pt]{article}
\usepackage{amsthm}
\usepackage{amsfonts}
\usepackage[T1]{fontenc}
\usepackage[latin1]{inputenc}
\usepackage{epsfig}

\usepackage[lflt]{floatflt}
\usepackage{epsfig}
\usepackage{graphicx}
\usepackage{amsmath}
\usepackage{amssymb}
\usepackage{color}
\usepackage{subfigure}
\usepackage[english]{babel}

\usepackage{makeidx}
\usepackage{latexsym}

\makeindex \textwidth=17cm \textheight=25cm \oddsidemargin=-0.5cm
\newtheorem{theorem}{Theorem}
\newtheorem{definition}{Definition}
\newtheorem{proposition}{Proposition}

\newtheorem{corollary}{Corollary}
\newtheorem{remark}{Remark}
\newcommand{\bbR}{{\mathbb R}}
\newcommand{\bbN}{{\mathbb N}}
\newcommand{\bbC}{{\mathbb C}}

\newcommand{\e}{\epsilon}
\newcommand{\al}{\alpha}
\newcommand{\la}{\lambda}
\newcommand{\p}{\partial}
\newcommand{\be}{\beta}
\newcommand{\G}{\Gamma}
\def\qed{\ \ \ \hbox{}\nolinebreak\hfill $ \Box     \  \  \  \  $ \par{}\medskip}

\begin{document}
\begin{center}\emph{}
\LARGE
\textbf{A Generalization of the Hopf's Lemma for the 1-D
Moving-Boundary Problem for the Fractional Diffusion
Equation and its Application to a Fractional Free-Boundary Problem.}
\end{center}
\medskip
\begin{center}
\normalsize
  Sabrina Roscani\\
\medskip
\small
Departamento de Matem\'{a}tica,  FCEIA, Universidad Nacional de Rosario, Pellegrini 250, Rosario, Argentina \\
CONICET, Argentina.\\
\textcolor{blue}{sabrina@fceia.unr.edu.ar} \\
\end{center}
\small

\small

\noindent \textbf{Abstract: }This paper deals with a theoretical mathematical analysis of a one-dimensional-moving-boundary problem for the time-fractional diffusion equation, where the time-fractional derivative of order $\al$ $\in (0,1)$ is taken in the Caputo's sense. A generalization of the Hopf's lemma is proved, and then this result is used to prove a monotonicity property for the free-boundary when a fractional free-boundary  Stefan problem is considered. \\


\noindent \textbf{MSC[2010] Primary:} 26A33,35R37,35R35; Secondary: 34K37, 35R11, 80A22.

\normalsize

\section{Introduction}

The development of the fractional calculus dates from the XIX century. Mathematicians as  Lacroix, Abel, Liouville, Riemann and Letnikov attemped to establish a definition of fractional derivative. But the definition given by Caputo in 1967, was an open door to the beginning of the physics applications, while the previous definitions  enabled a great theoretical development.\\

The research on the theory of fractional differential equations has begun to develop recently, and in the past decades many authors pointed out that derivatives and integrals of non-integer order are very useful in describing the properties of various real-world materials such as polymers or some types of non-homogeneous solids. The trend indicates that the new fractional order models are more suitable than integer order models previously used, since fractional derivatives give us an excellent tool for describing properties of memory and heritage of various materials and processes. Works in this direction are \cite{AST-libro, CuGi, GKS, Io2, SiUch}.

 This paper deals with the fractional diffusion equation (here in after FDE), obtained from the standard diffusion equation by replacing the first order time-derivative by a
fractional derivative of order $\alpha > 0 $ in the Caputo's sense:
$$   _0 D^{\alpha}_t u(x,t)=\lambda^2\,  u_{xx}(x,t), \quad  -\infty<x<\infty, \ t>0, \ 0<\al<1, $$
 where the fractional derivative in the Caputo's sense of arbitrary order $\al >0$ is given by
$$\,_{a} D^{\alpha}f(t)=\left\{\begin{array}{lc} \frac{1}{\Gamma(n-\al)}\int^{t}_{a}(t-\tau)^{n-\al-1} f^{(n)}(\tau)d\tau, &  n-1<\al<n\\
f^{(n)}(t), &   \al=n. \end{array}\right.$$
where  $n \in \bbN$ and $\G$ is the Gamma function defined by $\G(x)=\int_0^\infty  w^{x-1}e^{-w}dw$.\\

Whereas that 	the one-dimensional heat equation has become the paradigm for the all-embracing study of parabolic partial differential equations, linear or nonlinear (see Cannon \cite{Cannon}), the FDE plays a similar role in fractional parabolic operators.
\smallskip

 The FDE has been treated by a number of authors (see  \cite{FM-AnaPropAndAplOfTheW-Func,  Kilbas, Luchko3,FM-TheFundamentalSolution, Podlubny}) and, among the several applications that have been studied, Mainardi \cite{FM-libro} studied the application to the theory of linear viscoelasticity.\\
 
 Generalizations of the maximum principle for initial-boundary-value problems associated to the time-fractional diffusion equations were given by Luchko in \cite{Luchko3} and \cite{Luchko1}, and uniqueness results there were obtained.\\

Eberhard Frederich Ferdinand Hopf was an Austrian mathematician who made significant contributions in differential equations, topology and ergodic theory. One of his most important works are related to the strong maximum principle for partial differential equations of elliptic type (\cite{Hopf1}).\\

In his work \cite{Hopf2}, an important theorem related to the sign of the outside directional derivative of a function that is a solution to an elliptic partial differential inequality is proved. This theorem was
proved later for partial differential operators of parabolic type by A. Friedman \cite{AF-1958} and R. Viborni \cite{RV-1957} separately.
A weak adaptation of this theorem can be founded in \cite{Cannon}, named  Hopf's Lemma, and my propose is to generalize it for the FDE.\\ 

That is, under certain conditions that will be enunciated later, if $u$ is a solution of the following problem

 \begin{equation}\label{(P)}\left\{\begin{array}{lll}
         (i)\quad\,\,   _0D^{\al}_t u(x,t)=\lambda^2\,  u_{xx}(x,t)  & s_1(t)<x<s_2(t), \, 0<t\leq T, \, 0<\al<1 \\

        (ii)\quad   u(s_1(t),t)=g(t) &   0<t\leq T\\

          (iii)\,\, \,\,   u(s_2(t),t)=h(t) &  0<t\leq T\\

         (iv) \,\,\,\,\,  u(x,0)=f(x) &   a\leq x\leq b                   \end{array}\right. 
\end{equation}
where $s_1$ and $s_2$ are given, then, if $u$ assumes its maximum in a boundary point, let us say  $(s_2(t_0),t_0)$, it results that $u_x(s_2(t_0),t_0)>0$.

\section{A Fractional Hopf's Lemma} 

Let us consider the moving-boundary problem for the FDE defined in (\ref{(P)}) where:\\

(H1) The curve  $s_1$  is given and it is an upper Lipschitz continuous function.\\

(H2) The curve  $s_2$  is given and it is a lower Lipschitz continuous function.\\

(H3) $s_1(0)=a$, $s_2(0)=b$,  where  $a\leq b$,  and condition  $(iv)$ of problem $(\ref{(P)})$  is not considered if  $a=b$. \\

(H4) $s_1(t)<s_2(t) \, \forall t \in (0,T]$.\\

(H5) $f$  is a non-negative continuous function in  $[a,b]$.\\

(H6) $g$  and  $h$ are non-negative continuous functions in $(0,T]$.\\

\hspace{1 cm}

\noindent We will consider the following two regions:\\
 $D_T=\{ (x,t) / s_1(t)<x<s _2(t), \, 0<t\leq T \}$ and the so-called parabolic boundary  
 $\p_p D_T=\{ (s_1(t),t), 0<t\leq T \} \cup \{ (s_2(t),t), 0<t \leq T \} \cup \{(x,0), a\leq x \leq b\}$.

\begin{definition} 
  A function $u$ is a solution of problem 
  $(\ref{(P)})$  if $u=u(x,t)$ verifies the conditions in $(\ref{(P)})$ and 
\begin{enumerate}
	\item $u$ is defined in $[a_0,b_0]\times [0,T]$, where $a_0=\min\{ s_1(t), \, t \in [0,T]\}$ and $b_0=\max\{ s_2(t), \, t \in [0,T]\}$,
	\item $u\in $ $CW_{D_T}=C(D_T)\cap W^1((0,T))\cap C^2_x(D_T)$, 	where $W^1((0,T))=\{f\in C  ^1((0,T]) \, :\, f' \in L^1(0,T)\}$
\item $u$ is continuous in $D_T \cup \p_p D_T$ except perhaps at $(a,0)$ and $(b,0)$ where we will ask that 
 $$ 0\leq \underset{(x,t)\rightarrow (a,0)}{\liminf}u(x,t)\leq \underset{(x,t)\rightarrow (a,0)}{\limsup } u(x,t)<+\infty$$
  and 
	$$0 \leq \underset{(x,t)\rightarrow (b,0)}{\liminf}u(x,t)\leq \underset{(x,t)\rightarrow (b,0)}{\limsup } u(x,t)<+\infty. $$ 
\end{enumerate}

\end{definition}

\begin{remark} We ask $u$ to be defined in  $[a_0,b_0]\times [0,T]$ because the fractional derivative $_0D^\al_t u(x,t)$  involves the values of $u_t(x,\tau) $ for all $\tau$ in $[0,t]$. (See Figure 1).
\end{remark}

\begin{figure}[H]
\includegraphics[width=7cm]{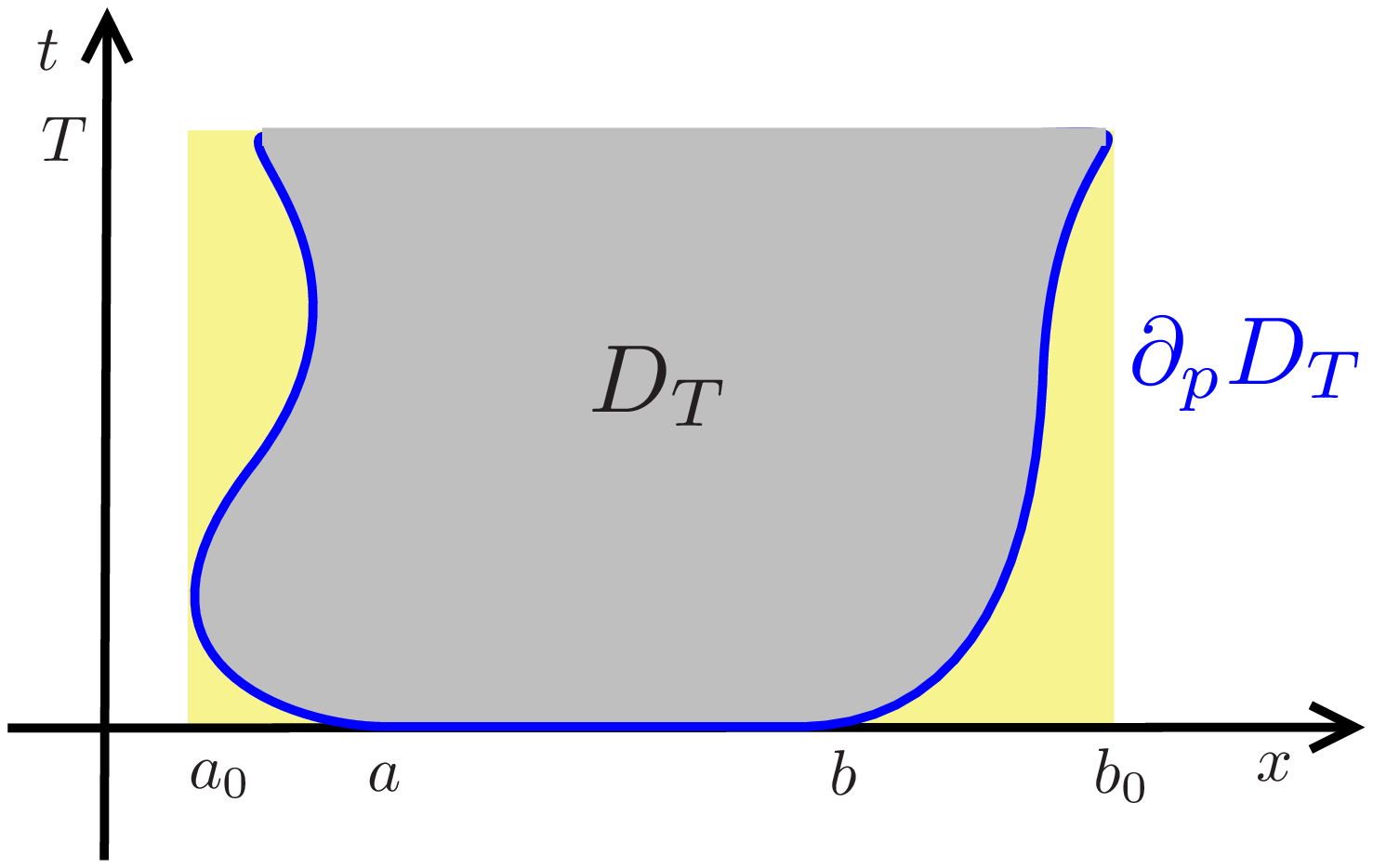}
\caption{Figure 1}
\end{figure}

\begin{remark} 
This kind of problems has not been studied in depth yet, but taking into account the results obtained in $\cite{RoSa1}$ and $\cite{RoSa2}$, where some fractional Stefan problems has been solved explicitly, it is easy to check that the following problem 

\begin{equation}\left\{\begin{array}{lll}
          _0D^{\al}_t u(x,t)=u_{xx}(x,t)  & 0<x<t^{\al/2}, \, 0<t\leq T, \, 0<\al<1 \\

          u(0,t)=B &   0<t\leq T\\

          u(t^{\al/2},t)=C &  0<t\leq T\\

                            \end{array}\right. 
\end{equation}

\noindent admits the solution given by

$$  u(x,t)=B+\frac{C-B}{1-W\left(-1,-\frac{\al}{2},1\right)}
           \left[1-W\left(-\frac{x}{t^{\al/2}},-\frac{\al}{2},1\right)\right],$$
where $W\left(\, \cdot\, ,-\frac{\al}{2},1\right)$ is the Wright function of parameters $\rho=-\frac{\al}{2}$
 and $\beta=1$ defined by

\begin{equation}\label{W} W(z;\rho;\be)=\sum^{\infty}_{k=0}\frac{z^{k}}{k!\G(\rho k+\be)}, \quad z\in \mathbb{C}, \, \rho>-1, \be\in \bbR.
\end{equation}

The function $1-W(-\, \cdot,-\frac{\al}{2},1)$ is the ``fractional error function'', so named because 
$$ \lim\limits_{\al\nearrow 1}1-W(-x,-\frac{\al}{2},1)=\, erf\left(\frac{x}{2}\right) $$
(see \cite{RoSa1} Theorem 4.1).
\end{remark}
\bigskip

\noindent Hereinafter I will take $\lambda=1$, I will call $D^\al$ to the fractional derivative in the Caputo's sense of extreme $a=0$, $_0D^\al_t$, and  $L^\al$ to the operator associated to the  FDE
\begin{equation}\label{L^al}
  L^\al:= \frac{\partial ^2}{(\partial x)^2 }- D^\al.
\end{equation}	

\bigskip

\begin{proposition} 
If $u$ is a function such that  $L^\al [u]>0$ in $D_T$, then $u$ can not attain its maximum at $D_T$.
\end{proposition}

\textit{Proof.}
Let us suppose that there exists $(x_0,t_0) \in D_T$ (that is $s_1(t_0)<x_0<s_2(t_0)$,  $0<t_0\leq T$), such that $u$ attains its maximum at $(x_0,t_0)$. \\
Due to the extremum principle for the Caputo derivative ( see  \cite{Luchko2}), we have that  
$ D^\al_t u(x_0,t_0)\geq 0 $. \\
On the other hand, since $ u \in C_x^2(D_T)$, $\frac{\partial^2 u }{\partial x^2}(x_0,t_0)\leq 0$. Then
$ L^\al[u](x_0,t_0)\leq 0 $, which is a contradiction.
\qed

\begin{corollary} \label{L[u]<0}
If $u$ is a function such that $L^\al [u]<0$ in $D_T$, then $u$ can not attain its minimum in $D_T$.
\end{corollary}

It is easy to adapt to the moving-boundary problem $(\ref{(P)})$ results obtained in \cite{Luchko1} for initial-boundary-value problems associated to the generalize FDE. For this reason we omit the proof of the following  result. 

\begin{theorem}\label{Ppio del maximo debil-negativo}
Let $u\in CW_{D_T}$ be a solution of $(\ref{(P)})$. Then either
$$  u(x,t)\geq  0 \, \forall \, (x,t) \in \overline{D_T}\quad \text{or  }  u \text{ attains its negative minimum on }  \p_p D_T. \,$$
\end{theorem}

\noindent Let us enunciate the main result of this paper.
\medskip

\begin{theorem}\label{frac hopf lemma s2} Let $u\in CW_{D_T}$ be a solution of problem $(\ref{(P)})$ satisfying the hypotheses $(H1)-(H6)$. 
 \begin{enumerate}
\item If there exists $t_0>0$  such that
\begin{equation}\label{(2-1)} u(s_2(t_0),t_0)=M=\sup_{\p_p D_T}{u},
\end{equation}                                                  
and,
\begin{equation}\label{(2-3)}\text{exists } \delta>0  \text{ such that }\, \left|s_1(t_0)-s_2(t_0)\right|\geq \delta \, \text{and } \, u(x,t_0)<M \, \forall x \in \, (s_2(t_0)-\delta, s_2(t_0)),  \end{equation}
then 
\begin{equation}\label{liminf u>0} \liminf_{x\nearrow s_2(t_0)}\frac{u(x,t_0)-u(s_2(t_0),t_0)}{x-s_2(t_0)} >0 .\end{equation}

If $u_x$ exists at $(s_2(t_0),t_0)$, then

\begin{equation}\label{u_x>0}u_x(s_2(t_0),t_0) >0  .\end{equation}

\item If there exists $t_0>0$  such that
\begin{equation}\label{(2-1)} u(s_2(t_0),t_0)=m=\inf_{\p_p D_T}{u},
\end{equation}                                                  
and,
\begin{equation}\text{exists } \delta>0 \text{ such that } \left|s_1(t_0)-s_2(t_0)\right|\geq \delta \, \text{and } u(x,t_0)>m \, \forall x \in \, (s_2(t_0)-\delta, s_2(t_0)),  \end{equation}
then 
\begin{equation}\label{limsup u<0} \limsup_{x\nearrow s_2(t_0)}\frac{u(x,t_0)-u(s_2(t_0),t_0)}{x-s_2(t_0)} <0 .\end{equation}

If $u_x$ exists at $(s_2(t_0),t_0)$, then

\begin{equation}\label{u_x<0}u_x(s_2(t_0),t_0) <0 . \end{equation}
\end{enumerate}
\end{theorem}
\textit{Proof.}
 I will prove 1. The proof of 2 is analogous. \\

Let us consider the function 
\begin{equation}\label{w_al} w_\al(x,t)=\e \left[1-\exp\{-\mu( x-s_2(t_0))\}\frac{E_\al(\mu A t^\al)}{E_{\al}(\mu A t_0^{\al})}\right]+M \end{equation}
where $A$, $\mu$ and $\e$ will be determined and $E_\al$ is the Mittag-Leffler function defined by 
\begin{equation}\label{E_al}
E_\al(z)=\sum_{k=0}^{\infty}\frac{z^k}{\G(\al k +1)}, \quad z\in \bbC, \al >0 .
\end{equation}

Note that $w_\al=M$ over the curve
\begin{equation}\label{curve} \exp\{-\mu( x-s_2(t_0))\}\frac{E_\al(\mu A t^\al)}{E_{\al}(\mu A t_0^{\al})}=1 .\end{equation}

Observe that the curve (\ref{curve}) is the graphic of the function  
\begin{equation} f(t)=\frac{1}{\mu}\ln\left(\frac{E_\al(\mu A t^\al)}{E_{\al}(\mu A t^{\al}_0)}\right)+s_2(t_0), \quad t\in (0,t_0] \end{equation}
Clearly $f(t_0)=s_2(t_0)$ and it is easy to check that 
\begin{equation}\label{f crec}f \text{ is an increasing function if } \mu >0. \end{equation}

Our next goal is to prove that there exists  $t_1<t_0$ such that $f(t)< s_2(t)$ $\forall t\in (t_1,t_0)$. \\

Due to  $(H2)$, $s_2$ is a lower Lipschitz continuous function, then there exists a constant $L>0$  such that  
$$ \frac{s_2(t)-s_2(t_0)}{t-t_0}\leq L \quad \forall \,\, 0\leq t\leq t_0. $$ Therefore
\begin{equation}\label{2-5}  s_2(t)\geq L(t-t_0)+s_2(t_0) \quad \forall \,\, 0\leq t\leq t_0. \end{equation}

Taking into account that  $E_{\al}(\mu A t^{\al})=\sum_{k=0}^{\infty}\frac{\left(\mu A t^{\al}\right)^k}{\G(\al k+1)}$ is an uniform convergent series over compact sets contained in $(0,t_0]$ and that  $z\G(z)=\G(z+1)$ $\forall \, z\in \bbC $, we have that
$$ \left[E_{\al}(\mu A t^{\al})\right]'=\sum_{k=1}^{\infty}\frac{(\mu A)^k\al k t^{\al k -1}}{\G(\al k+1)}= \sum_{k=0}^{\infty}\frac{(\mu A)^{k+1} t^{\al k +\al-1}}{\G(\al k+\al)}=\mu A t^{\al-1}E_{\al,\al}(\mu A t^{\al}),$$
where the function $E_{\al,\al}$ is the generalized Mittag-Leffler function of parameters $\rho=\beta=\al$ defined by
\begin{equation}\label{E_rho,be}
E_{\rho,\be}(z)=\sum_{k=0}^{\infty}\frac{z^k}{\G(\rho k +\be)}, \quad \,z \in \bbC, \rho>0, \beta \in \bbC.
\end{equation}

Then, 

\begin{equation}\label{2-6}f'(t)=\frac{1}{\mu}\frac{1}{E_{\al}(\mu A t^{\al})}\frac{\mu A}{t^{1-\al}}E_{\al,\al}(\mu A t^{\al})=\frac{A}{t^{1-\al}}\frac{E_{\al,\al}(\mu A t^{\al})}{E_{\al}(\mu A t^{\al})}.
\end{equation}

Now, let us define the function $H:\bbR^+_0\rightarrow \bbR /$ $H(t)=\frac{E_{\al,\al}(\mu A t^{\al})}{E_{\al}(\mu A t^{\al})}$. \\

$H$  is a positive function and it is a quotient of continuous functions, where its denominator is grater than $1$ in $[0,\infty)$, then $H$ is continuous in  $[0,\infty)$.\\
$H(0)=\frac{1}{\G(\al)}>0$ because $0<\al<1$.\\
$H(+\infty)= C>0 $ because it is a quotient of continuous functions with equal order in $\infty$ (see \cite{FM-Mittag-Leffler}).\\
Then, we can assure that there exists  $m_0>0 $ such that 
\begin{equation}\label{2-7} H(t)\geq m_0 \quad \forall \, t\geq 0, \forall \, A, \mu >0.\end{equation}

From  (\ref{2-6}) and (\ref{2-7}), we have that

\begin{equation}\label{cota f'} 
f'(t_0)\geq \frac{A}{t_0^{1-\al}}m_0 .
\end{equation}

Selcting  $A>0$ such that  $\frac{A}{t_0^{1-\al}}m_0>L$ we can assure that 

\begin{equation}\label{2-8} f'(t_0)>L.\end{equation}

Lately, let be $\rho >0$ / $f'(t_0)-\rho >L$. Due to the differentiability of $f$ at $t_0$ we can assure that there exists  $t_1<t_0$ such that $\forall \, t \, \in (t_1,t_0)$, 
$$ L<f'(t_0)-\rho<\frac{f(t)-f(t_0)}{t-t_0}\Rightarrow f(t)<L(t-t_0)+f(t_0)=L(t-t_0)+s_2(t_0)\leq s_2(t) $$
  Note that due to $(\ref{(2-3)})$ and (\ref{f crec}), we can select  $t_1$ so that  $s_1(t)<f(t)< s_2(t)$ $\forall t \in (t_1,t_0)$.\\

 Now, let be $A(x_1,t_0)$ (where $x_1=f(t_1)$), $B(s_2(t_0),t_0)$ 
and $C(x_1,t_1)$. Hypothesis $(\ref{(2-3)})$ allows us to set $t_1$ again such that $x_1 \in \, (s_2(t_0)-\delta,s_2(t_0))$  and 
$u<M$ in $\overline{AC}$. \\
 Let be $\mathcal{R}$  the region limited by $\overline{AB}$, $\overline{AC}$  and the portion of graph of $f$ from  $B$ to $C$, which we will  call $\widehat{CB}$.(See Figure 2)\\

The region $\mathcal{R}_{t_0}=\mathcal{R}^\circ \cup (\overline{AB}-\left\{A,B\right\})$ and its parabolic boundary
$\p_p\mathcal{R}=\overline{AC}\cup \widehat{CB}$ will be considered.\\

\begin{figure}[H]
 \includegraphics[width=7cm]{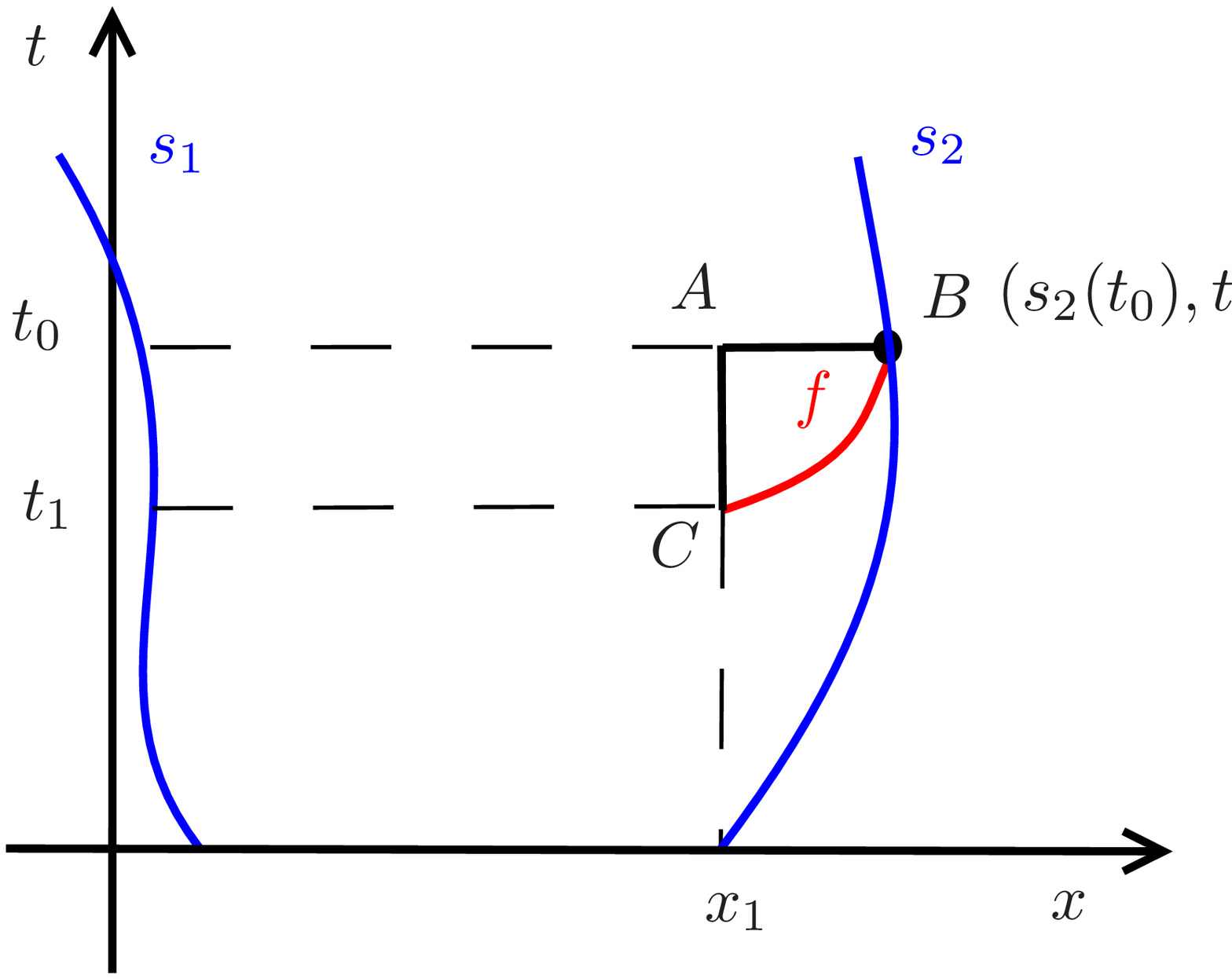}
\caption{Figure 2}
\end{figure}  

 Next, we will analyze the behavior of  $u$ and $w_\al$ in the parabolic boundary  $\p_p\mathcal{R}$.\\
 Let be $M_0=\underset{t_1\leq t\leq t_0}{max}{u(x_1,t)}$. Because of the continuity of $u$, the hypothesis $(\ref{(2-3)})$ and resetting $t_1$ if it is necessary, we can affirm that $M_0<M$. Calling $\eta=M-M_0$, it yields that 
\begin{equation}\label{2-8} u\leq M-\eta \quad \text{ in } \overline{AC}.\end{equation}

\begin{equation}\label{2-9} u\leq M \quad \text{ in } \widehat{CB}.\end{equation}
 By the other side, 
\begin{equation}\label{2-9'} w_\al =M \quad \text{ in }  \widehat{CB}.\end{equation}
 In $\overline{AC}$, considering that $E_\al(\mu A t^\al)$  is an increasing function, we have that, 
\begin{equation}\label{2-10'} w_\al (x_1,t)=\e \left[1-\exp\left\{ -\mu (x_1-s_2(t_0))\right\}\frac{E_\al(\mu A t^\al)}{E_{\al}(\mu A t_0^{\al})} \right]+M\geq \e 
\left[1-\exp\left\{ -\mu (x_1-s_2(t_0))\right\}\right]+M.\end{equation}
 
Taking $\e=\frac{\eta}{ exp\left\{ -\mu (x_1-s_2(t_0))\right\}-1}$, it results that 
\begin{equation}\label{2-10} w_\al (x_1,t)\geq -\eta +M. \end{equation}

From \cite{Kilbas} we have that 
\begin{equation}\label{D_al E}
 D^\al (E_\al(\mu A t^\al))=\mu A E_\al(\mu A t^\al).
\end{equation}

Now, applying the operator $L^\al$ to the function $w_\al$ and using $(\ref{D_al E})$ it yields that  

\begin{equation}\label{2-11} L^\al[w_\al](x,t)=\e \exp\left\{ -\mu (x-s_2(t_0))\right\} \frac{E_\al(\mu A t^\al)}{E_{\al}(\mu A t_0^{\al})} (\mu A-\mu^2)<0 \quad \text{ if } \mu=A+1. \end{equation}

Finally, we define the function $z=w_\al -u$ in $\mathcal{R}$. Let us analyze the behavior of $z$ in the parabolic boundary $\p_p\mathcal{R}$. \\

From (\ref{2-8}) and (\ref{2-10}), $z\geq 0$ in $\overline{AC}$.\\

From (\ref{2-9}) and (\ref{2-9'}), $z\geq 0$ in $\widehat{CB}$ .\\

Also, $L^{\al}[z]=L^{\al}[w_\al]-L^\al[u]<0 $ in $\mathcal{R}_{t_0}$. \\

Applying Corollary \ref{L[u]<0}, we can state that  $z$ cannot assume its minimum at $\mathcal{R}_{t_0}$. Then  
$$ z\geq 0 \quad \text{ in } \, \mathcal{R}$$

In particular, 
\begin{equation}\label{(2-12)} z(x,t_0)=w_\al(x,t_0)-u(x,t_0)\geq 0, \quad \forall \,x_1 \leq x\leq s_2(t_0). \end{equation}

Recalling that $u(s_2(t_0),t_0)=w_\al(s_2(t_0,t_0))=M$, the next expression es equivalent to (\ref{(2-12)}):

\begin{equation}\label{(2-13)}\frac{u(x,t_0)-u(s_2(t_0),t_0)}{x-s_2(t_0)}\geq \frac{w_\al(x,t_0)-w_\al(s_2(t_0,t_0))}{x-s_2(t_0)}.\end{equation}

Then
$$ \liminf_{x\nearrow s_2(t_0)}\frac{u(x,t_0)-u(s_2(t_0),t_0)}{x-s_2(t_0)} \geq \liminf_{x\nearrow s_2(t_0)}\frac{w_\al(x,t_0)-w_\al(s_2(t_0,t_0))}{x-s_2(t_0)}. $$

But $w_\al$ is a differentiable function at $(s(t_0),t_0)$, then 
$$ \liminf_{x\nearrow s_2(t_0)}\frac{w_\al(x,t_0)-w_\al(s_2(t_0,t_0))}{x-s_2(t_0)}=(w_\al)_x(s_2(t_0),t_0)=\e \mu=\e (A+1)>0 $$
and $(\ref{liminf u>0})$ holds.\\

Finally, if $u_x$ exists at $(s_2(t_0),t_0)$, $(\ref{(2-13)})$ implies that 

\begin{equation}\label{2-14}
u_x(s_2(t_0),t_0) \geq (w_\al)_x(s_2(t_0),t_0)>0
\end{equation}
  
and $(\ref{u_x>0})$ holds.

\qed

The same result is valid if we consider $s_1$ instead of $s_2$.

\begin{theorem} Let $u\in CW_{D_T}$ be a solution of problem $(\ref{(P)})$ satisfying the hypotheses $(H1)-(H6)$.  
 \begin{enumerate}
\item If there exists $t_0>0$  such that
\begin{equation}\label{(2-1)} 
u(s_1(t_0),t_0)=M=\sup_{\p_p D_T}{u}, \end{equation}                                      
and,
\begin{equation}\text{exists } \delta>0 \text{ such that }\left|s_1(t_0)-s_2(t_0)\right|\geq \delta \text{and } u(x,t_0)<M \,\forall x \in (s_1(t_0), s_1(t_0)+\delta),  \end{equation}
then 
\begin{equation} \limsup_{x\nearrow s_1(t_0)}\frac{u(x,t_0)-u(s_1(t_0),t_0)}{x-s_1(t_0)} <0 .\end{equation}
If $u_x$ exists at $(s_1(t_0),t_0)$, then
\begin{equation}u_x(s_1(t_0),t_0) <0  .\end{equation}

\item If there exists $t_0>0$  such that
\begin{equation}\label{(2-1)} u(s_1(t_0),t_0)=m=\inf_{\p_p D_T}{u},
\end{equation}                                                
and,

\begin{equation}\text{exists } \delta>0 \, \text{ such that } \left|s_1(t_0)-s_2(t_0)\right|\geq \delta \, \text{and } u(x,t_0)>m \, \forall x  \in  (s_1(t_0), s_1(t_0)+\delta),   \end{equation}
then 
\begin{equation} \liminf_{x\nearrow s_1(t_0)}\frac{u(x,t_0)-u(s_1(t_0),t_0)}{x-s_1(t_0)} >0 .\end{equation}

If $u_x$ exists at $(s_1(t_0),t_0)$, then

\begin{equation}u_x(s_1(t_0),t_0) >0 . \end{equation}
\end{enumerate}
\end{theorem}

\section{An Application to Fractional Free-Boundary Stefan Problems. }

Let us consider now the following fractional free-boundary Stefan problem for the FDE, where we have replaced the Stefan condition $\dfrac{d s(t)}{dt}=k u_x(s(t),t), \, t>0,$ by the fractional Stefan condition
$$   D^{\alpha}s(t)= -k u_x(s(t),t),\quad t>0, \quad 0<\al<1. $$
\begin{equation}{\label{St1}}
\left\{\begin{array}{lll}
     (i)  \quad \,\, \,    D^{\al} u(x,t)=u_{xx}(x,t) &   0<x<s(t), \,  0<t<T, \,  0<\al<1 , \, \, \la>0\\
     (ii) \,\, \quad  u(x,0)=f(x) & 0\leq x\leq b=s(0)\\ 
       (iii)  \quad  u(0,t)=g(t)   &  0<t\leq T  \\

         (iv) \, \quad u(s(t),t)=0 & 0<t\leq T\\

        (v) \quad \, \, D^{\al}s(t)=-k u_x(s(t),t) & 0<t\leq T, \quad  k>0 \, \text{ constant}
                                             \end{array}\right.\end{equation}

\begin{definition}\label{Def sol St} A pair $\{u,s\} $ is a solution of  problem $(\ref{St1})$  if

\begin{enumerate}
    \item $u$ and $s$ satisfies $(\ref{St1})$,
   	\item $s$ is a continuous function in $[0,T]$ such that $s \in W(0,T)=\{ f \in  C^1((0,T]) \, / \, f' \in L^1(0,T) \}. $
   	\item $u$ is defined in $[0,b_0]\times[0,T]$ where $b_0= \max\{s(t), 0 \leq t \leq T \}$,
	\item $u\in $ $CW_{D_T}$.
\item $u$ is continuous in $D_T \cup \p_p D_T$ except  perhaps at $(0,0)$ and $(b,0)$ where we will ask that 
 $$ 0\leq \underset{(x,t)\rightarrow (0,0)}{\liminf}u(x,t)\leq \underset{(x,t)\rightarrow (0,0)}{\limsup } u(x,t)<+\infty$$
  and 
	$$0\leq \underset{(x,t)\rightarrow (b,0)}{\liminf}u(x,t)\leq \underset{(x,t)\rightarrow (b,0)}{\limsup } u(x,t)<+\infty. $$ 
    \item There exists $u_x(s(t),t)$ for all $t \in (0,T]$.
    \end{enumerate}
\end{definition}

This kind of problems have been recently treated in \cite{Atkinson, Garra,Li-Xu-Jiang, RoSa1}), and our goal now is to prove the next theorem involving the monotonicity of the free boundary.

\begin{theorem}
Let  $\{u_1,s_1\}$ and $\{u_2,s_2\}$ be solutions of the fractional free-boundary Stefan  problems $(\ref{St1})$ corresponding, respectively, to the data $\{b_1, f_1, g_1\}$ and $\{ b_2,f_2, g_2 \}$.
Suppose that $ b_1<b_2$, $0\leq f_1\leq f_2$ and $0\leq g_1\leq g_2$. Then  $s_1(t)<s_2(t) \forall \, t\in [0,T)$.

\end{theorem}

\textit{Proof.}
$s_1$ and $s_2$ are continuous functions and $s_1(0)=b_1<b_2=s_2(0)$. Suppose that the set $A=\left\{t\in [0,T]/\, (s_1-s_2)(t)=0\right\}\neq \emptyset$, and let be $t_0=\min{A}$. Due to the continuity of $s_1$ and $s_2$, $s_1(t_0)=s_2(t_0)$, and $t_0$ is the first $t$ for which $s_1(t_0)=s_2(t_0)$. \\
 Let be $h:[0,t_0]\rightarrow \bbR$ the function $h(t)=(s_1-s_2)(t)$. $h$ has the following properties:\\
($h$-1) \hspace{0.5 cm} $h \in C^1(0,t_0]\cap C[0,t_0]$ (due to definition \ref{Def sol St}). \\
($h$-2) \hspace{0.5 cm} $h(0)=b_1-b_2<0$. \\
($h$-3) \hspace{0.5 cm} $h$ is a non positive function and $h(t_0)=0$.\\
From ($h$-1)-($h$-3), $h$ attains its maximum value at $t_0$.\\
 Using the estimate ( eq. (12) of Theorem 2.1 of \cite{AlRefai-Luchko}), it results that
\begin{equation}\label{acotacion D^al}
D^\al h (t_0)\geq \frac{h(t_0)-h(0)}{t_0^\al \G(1-\al)}.
\end{equation}
Then,
\begin{equation}\label{Dh>0}
D^\al h (t_0)\geq \frac{b_2-b_1}{t_0^\al \G(1-\al)}>0.
\end{equation}
    
Taking into account the linearity of the Caputo's fractional derivative and  that $s_1$  and $s_2$ satisfies the Stefan condition (\ref{St1})-(v), (\ref{Dh>0}) implies that 
\begin{equation}\label{mono-2}
u_{2x}(s_2(t_0),t_0)- u_{1x}(s_1(t_0),t_0)>0.
\end{equation}

By the other hand, the function $w(x,t)=u_2(x,t)-u_1(x,t)$ is a solution of the following moving-boundary problem
\begin{equation}{\label{mov-1}}
\left\{\begin{array}{lll}
          D^{\al} w(x,t)=w_{xx}(x,t) &   0<x<s_1(t), \,  0<t\leq t_0, \,  0<\al<1 , \, \\
  
          w(0,t)=(g_2-g_1)(t)\geq 0   &  0<t\leq t_0  \\

          w(s_1(t),t)=u_2(s_1(t),t) & 0<t\leq t_0\\
                           
   w(x,0)=(f_2-f_1)(x)\geq 0 & 0\leq x\leq b_1=s_1(0)                        \end{array}\right.\end{equation}
 Applying Theorem \ref{Ppio del maximo debil-negativo} to $u_2$ in the region $\overline{D^2_{t_0}}$ where $D^2_{t_0}=\{ (x,t) \, /\, 0<t\leq t_0, \,0<x<s_2(t) \}$, it results  that $u_2(s_1(t),t)\geq 0.$ \\
 Now, applying again Theorem  \ref{Ppio del maximo debil-negativo} to $w$ in problem (\ref{mov-1}), it results that  $w\geq 0$ in $\overline{D^1_{t_0}}$, where $D^1_{t_0}=\{ (x,t) \, /\, 0<t\leq t_0, \,0<x<s_1(t) \}$. \\
 Then we can state that $w$ attains a minimum at $(s_1(t_0),t_0)$.\\
 Now, if there exists $\e >0$ such that $w(x,t_0)>0$ for all  $x\in (s_1(t_0)-\e,s_1(t_0))$, applying Theorem  \ref{frac hopf lemma s2}-2 we can conclude that $w_x(s_1(t_0),0)<0$. And then
\begin{equation}\label{mono-3}
 u_{2x}(s_2(t_0),t_0)- u_{1x}(s_1(t_0),t_0)<0, 
\end{equation} 
which contradicts (\ref{mov-1}).\\
If, by contrast, we have sequence $\{\e_n \}$ such that $\e_n\rightarrow 0$ and for every $n\in \bbN$ there exists $x_n \in (s_1(t_0)-\e_n,s_1(t_0))$ such that $w(x_n,t_0)=0$, then 
$$\lim_{n\rightarrow  \infty} \frac{w(x_n,t_0)-w(s_1(t_0),t_0)}{\e_n}=0 $$
and due to the existence of $w_x(s_1(t_0),t_0)$ (definition  \ref{Def sol St}-6), it results that $w_x(s_1(t_0),t_0)=0 $. Then

\begin{equation}\label{mono-4}
 u_{2x}(s_2(t_0),t_0)- u_{1x}(s_1(t_0),t_0)=0, 
\end{equation} 
which contradicts (\ref{mov-1}) again.\\
This contradiction comes from assuming that there exists $t_0>0$ such that $t_0$ is the first $t$ for which $s_1(t_0)=s_2(t_0)$. Therefore $s_1(t)<s_2(t) \forall \, t\in [0,T)$. 
\qed

\section{Conclusions}
A one-dimensional moving-boundary problem for the time-fractional diffusion equation was presented, where the time-fractional derivative of order $\al$ $\in (0,1)$ was taken in the Caputo's sense. Then, a generalization of the Hopf's lemma, involving the behavior of the partial $x$-derivative of the solution at a boundary point, was  proved. This last result was used to prove a monotonicity property of the free-boundary, when a free-boundary Stefan problem was considered.


\section{Acknowledgments}

\noindent This paper has been sponsored by PIP N$^\circ$ 0534 from CONICET-UA  and ING349, from Universidad Nacional de Rosario, Argentina.


\bibliographystyle{plain}
\bibliography{Roscani_bibfile}

\end{document}